\documentclass[smallextended, amsart]{amsproc}
\usepackage{mathptmx}
\usepackage{amsmath,amssymb,enumerate,epic,graphicx,mathrsfs}
\usepackage{amssymb,amsthm,amsmath,latexsym}
\newtheorem{theorem}{\sc Theorem}[section]
\newtheorem{lemma}[theorem]{\sc Lemma}

\newtheorem{corollary}[theorem]{\sc Corollary}

\newtheorem{conjecture}[theorem]{\sc Question}

\title{A new characterization of simple $K_3$-groups using same-order type}
\author{Igor Lima and Josyane Pereira}
\address{Departamento de Matem\'atica, Universidade de Bras\'ilia,
Brasilia-DF, 70910-900 Brazil }
\email{(Lima) igor.matematico@gmail.com, (Pereira) josyanedsp@gmail.com}

\thanks{The first author was supported by DPI/UnB, Brazil.}
\thanks{The second author was supported by CNPq.}


\date{May 2021}

\usepackage{graphicx}

\begin{document}

\maketitle

\begin{abstract}
Let $G$ be a group, define an equivalence relation $\sim$ as below:
$$\forall \ g, h \in G, \ g \sim h \Longleftrightarrow |g| = |h|$$
the set of sizes of equivalence classes with respect to this relation is called the same-order type of $G$ and denoted by $\alpha(G)$. And $G$ is said a $\alpha_n$-group if $|\alpha(G)|=n$. Let $\pi(G)$ be the set of prime divisors of the order of $G$. A simple group of $G$ is called a simple $K_n$-group if $|\pi(G)|=n$. We give a new characterization of simple $K_3$-groups using same-order type. Indeed we prove that a nonabelian simple group $G$ has same-order type \{r, m, n, k, l\} if and only if $G \cong PSL(2,q)$, with $q=7, 8$ or $9$. This result generalizes the main results in \cite{KKA}, \cite{Sh} and \cite{TZ1}. Moreover based on the main result in \cite{TZ1} we have the natural question: \textit{Let $S$ be a nonabelian simple $\alpha_n$-group and $G$ a $\alpha_n$-group such that $|S|=|G|$. Then $S \cong G$}. In this paper with a counterexample we give a negative answer to this question.
    
\end{abstract}

\section{Introduction and Preliminaries}

In this paper all the groups we consider are finite. 

Let $G$ a group and $\pi_e(G)$ be the set of element orders of $G$. Let $t \in \pi_e(G)$ and $s_t$ be the number of elements of order $t$ in $G$. Let $nse(G) = \{ s_t | \ t \in \pi_e(G) \}$  the set of sizes of elements with the same order in $G$. Some authors have studied the influence of $nse(G)$ on the structure of $G$ (see \cite{An}, \cite{Mal}, \cite{TZ1} and \cite{Za}). For instance R. Shen in \cite{Sh} proved that $A_4 \cong PSL(2, 3)$, $A_5 \cong PSL(2, 4) \cong PSL(2, 5)$ and $A_6 \cong PSL(2, 9)$ are uniquely determined by $nse(G)$. As a continuation in \cite{KKA} was proved that if $G$ is a group such that $nse(G) = nse(PSL(2, q))$, where $q \in \{7, 8, 11, 13 \}$, then $G \cong PSL(2, q)$. In \cite{SSM} and \cite{TZ1} new characterizations of $A_5$ were given using $nse(A_5)$. The authors in \cite{SSM} proved that $A_5$ is the only group such that $nse(A_5) = \{ 1, 15,20,24 \}$ and the authors in \cite{TZ1} generalized that a nonabelian simple group $G$
has same-order type $\{r, m, n, k \}$ if and only if $G \cong A_5$ (see Th. 1.1 \cite{TZ1}). 

Let $G$ be a group, in \cite{TZ1} was defined an equivalence relation $\sim$ as below:
$$\forall \ g, h \in G, \ g \sim h \Longleftrightarrow |g| = |h|$$
the set of sizes of equivalence classes with respect to this relation is called the same-order type of $G$ and denoted by $\alpha(G)$. And $G$ is said a $\alpha_n$-group if $|\alpha(G)|=n$. Note that $\alpha(G)$ is equal to the set of sizes of elements with the same order in $G$, hence $|nse(G)|=|\alpha(G)|$.

We give a new characterization of $PSL(2,7), PSL(2,8)$ and $PSL(2,9)$ using same-order type.

\begin{theorem} Let $G$ be a simple $K_3$-group with same-order type \{r, m, n, k, l\}. Then  $G \cong PSL(2,7), PSL(2,8)$ or $PSL(2,9)$.

\label{main}
\end{theorem}

This result generalizes the main results in \cite{KKA}, \cite{Sh} and \cite{TZ1}. Combination the main results in \cite{KKA} and \cite{Sh} with Theorem \ref{main} we have the following result

\begin{corollary} A simple $K_3$-group $G$ has same-order type $\{r, m, n, k, l\}$ if and only if $G \cong PSL(2,7), PSL(2,8)$ or $PSL(2,9)$.

\label{coro}
\end{corollary}

We see easily that the only $\alpha_1$-groups are $1$ and a cyclic group of order $2$. In \cite{Sh} R. Shen characterized $\alpha_2$-groups as nilpotent groups and $\alpha_3$-groups as solvable groups. Moreover Taghvasani-Zarrin (see Th. 1.1 in \cite{TZ1}) showed that the only nonabelian simple $\alpha_4$-group is the $A_5$. As noted in \cite{KKA} and \cite{TZ1} finite groups $G$ cannot be determined by $nse(G)$. Indeed in 1987 Thompson gave a first example as follows: Let $G_1 = (C_2 \times C_2 \times C_2 \times C_2) \rtimes A_7$ and $G_2 = PSL(3,4) \rtimes C_2$ be the maximal subgroups of Mathieu group $M_{23}$. Then $nse(G_1) = nse(G_2)$, but $G_1 \ncong G_2$.

Motived by the main result in \cite{TZ1} about a new characterization of $A_5$ using same-order type, we have the natural question

\begin{conjecture} Let $S$ be a nonabelian simple $\alpha_n$-group and $G$ a $\alpha_n$-group such that $|S|=|G|$. Then $S \cong G$.

\label{conje}
\end{conjecture}

We give a negative answer to this question in the last section.

\section{Proof of Theorem \ref{main}}

We need of one preliminary result to prove the main Theorem. The following result is a property very interesting of simple groups (see Lemma 2.7 in \cite{TZ1}).

\begin{lemma} Let $G$ be a nonabelian simple group. Then there exist two odd prime divisors $p$ and $q$ of the order of $G$ such that $s_p \neq s_q$.

\label{oddprime}
\end{lemma}

In fact if $G$ is a nonabelian simple group then there exist two odd prime divisors $p$ and $q$ of the order of $G$ such that $\{ 1, s_2, s_p, s_q\} \subseteq \alpha(G)$ (see Corollary 2.8 in \cite{TZ1}).

We are now ready to conclude the proof of main Theorem.\\

\textbf{\proofname{ of Theorem \ref{main}:}} As $G$ is a nonabelian simple group, it follows that $s_2 > 1$, w.l.g. $r=1$ and $s_2=m$. From Lemma \ref{oddprime} there exist odd prime divisors $p$ and $q$ of the order of $G$ such that $n=s_p \neq s_q=k$, hence $\pi(G) = \{2,p,q \}$ because $G$ is a simple $K_3$-group. Therefore $\{1, s_2, s_p, s_q \} \subseteq \alpha(G) = \{r, m, n, k, l\}$. So there exist a divisor $t \notin \pi(G)$ of order of $G$ such that $s_t=l$. It's well known that the only nonabelian simple groups of order divisible by exactly three primes are the following eight groups: $PSL(2,q)$, where $q \in \{5, 7, 8, 9, 17\}, PSL(3,3), PSU(3,3)$ and $PSU(4,2)$, see Th. 1 and Th. 2 in \cite{Her}. Now we arguing as in the proof of Th. 1.1 in \cite{TZ1} and a GAP check yields that $|\alpha(PSL(2,7))|=5$, $|\alpha(PSL(2,8))|=5$, $|\alpha(PSL(2,9))|=5$ and all others groups are $\alpha_n$-groups with $n \geq 6$ (except $A_5$ since $|\alpha(A_5)|=4$). The result is follows.

\section{A counterexample to a Question \ref{conje}}

Now we give a counterexample to the Question \ref{conje}. Firstly we observed that by the main Theorem in \cite{Sh}, we have that $\alpha(PSL(2,7))$ is uniquely determined and we have that $\alpha(PSL(2,7)) = \{1,21,56,42,48\}$ hence $PSL(2,7)$ is a $\alpha_5$-group. Let $G=Q_8 \times (C_7 \rtimes C_3)$, where $Q_8$ is the quaternion group of order $8$. As $|G|=168$ and $G$ is a soluble group then is sufficient to prove that $|\alpha(G)|=5$. Indeed the only $2$-Sylow subgroup $Q_8$ is a normal subgroup of $G$ and using Sylow's Theorem it follows that $s_2=8$. Note that a $7$-Sylow subgroup of $G$ is isomorphic to $C_7$ and is a normal subgroup of $Q_8 \cdot C_7$ and $C_7 \rtimes C_3$, hence the normalizer $N$ of $C_7$ has same order of $G$. Again from Sylow's Theorem we have that $C_7$ is a normal subgroup of $G$ and $s_7=56$. As the number of $3$-Sylow subgroup of $G$ is $7$, then $s_3=14$. The number of elements of $G$ of order $2$, $4$ are respectively $1$ and $6$, hence $s_2=1$, $s_4=6$ and consequently  $s_6=14$, $s_{12}=84$, $s_{14}=6$ and $s_{28}=36$ (because of direct product in the structure of $G$).  Therefore $\alpha(G)= \{1,1,14,6,14,6,84,6,36 \}$ and $G$ is a $\alpha_5$-group. Clearly $|PSL(2,7)|=168=|G|$ but $PSL(2,7) \ncong G$.

We can obtain others groups $G$ with the computational group theory system GAP \cite{Gap}: $G= C_7 \times (Q_8 \rtimes C_3)$ or $G= C_2 \times ((C_{14} \times C_2) \rtimes C_3)$. These groups are also counterexamples to the Question \ref{conje}.
\\

{\bfseries Acknowledgments}: We thank L. J. Taghvasani for pointing out a mistake of ours and talking to us about Question \ref{conje} in the olderst version of \cite{TZ1}.

\end{document}